\documentclass[11pt,reqno]{amsart}

\usepackage{amssymb,latexsym}
\usepackage{graphicx}
\usepackage{url}		

\calclayout
\makeatletter
\newcommand{\monthyear}[1]{%
  \def\@monthyear{\uppercase{#1}}}
\newcommand{\volnumber}[1]{%
  \def\@volnumber{\uppercase{#1}}}
\AtBeginDocument{%
\def\ps@plain{\ps@empty
  \def\@oddfoot{\@monthyear \hfil \thepage}%
  \def\@evenfoot{\thepage \hfil \@volnumber}}
\def\ps@firstpage{\ps@plain}
\def\ps@headings{\ps@empty
  \def\@evenhead{%
    \setTrue{runhead}%
    \def\thanks{\protect\thanks@warning}%
    \uppercase{The Fibonacci Quarterly}\hfil}%
  \def\@oddhead{%
    \setTrue{runhead}%
    \def\thanks{\protect\thanks@warning}%
    \hfill\uppercase{Tribonacci squares identities}}%
  \let\@mkboth\markboth
  \def\@evenfoot{%
    \thepage \hfil \@volnumber}%
  \def\@oddfoot{%
    \@monthyear \hfil \thepage}%
  }%
\footskip=25pt
\pagestyle{headings}%
}
\makeatother

\theoremstyle{plain}
\numberwithin{equation}{section}
\newtheorem{thm}{Theorem}[section]

\newtheorem{lemma}[thm]{Lemma}

\newtheorem{idn}[thm]{Identity}

\begin{document}
\bibliographystyle{fq} 

\monthyear{}
\volnumber{Volume ??, Number ?}
\setcounter{page}{1}

\title{Identities involving the tribonacci numbers squared via tiling
  with combs}
\author{Michael A. Allen*} 
\author{Kenneth Edwards}
\address{Physics Department, Faculty of Science,
Mahidol University, Rama 6 Road, Bangkok 10400 Thailand}

\email{maa5652@gmail.com \textrm{(*corresponding author)}}
\email{kenneth.edw@mahidol.ac.th}

\begin{abstract}
The number of ways to tile an $n$-board (an $n\times1$ rectangular
board) with $(\frac12,\frac12;1)$-, $(\frac12,\frac12;2)$-, and
$(\frac12,\frac12;3)$-combs is $T_{n+2}^2$ where $T_n$ is the $n$th
tribonacci number.  A $(\frac12,\frac12;m)$-comb is a tile composed of
$m$ sub-tiles of dimensions $\frac12\times1$ (with the shorter sides
always horizontal) separated by gaps of dimensions $\frac12\times1$.
We use such tilings to obtain quick combinatorial proofs of identities
relating the tribonacci numbers squared to one another, to other
combinations of tribonacci numbers, and to the Fibonacci, Narayana's
cows, and Padovan numbers. Most of these identities appear to be new.
\end{abstract}

\maketitle

\section{Introduction}
It is well-known that a combinatorial interpretation of the Fibonacci
number $F_{n+1}$ (where $F_n=F_{n-1}+F_{n-2}$, $F_1=1$, $F_0=0$) is
the number of ways to tile an $n$-board (a linear array of $n$ unit
square cells) with tiles of size $1\times1$ and $1\times2$
\cite{BQ=03}. As a generalization of this, combs were recently
introduced to give a combinatorial interpretation of
$s_n^rs_{n+1}^{p-r}$ where $s_n=\sum_{i=1}^q v_i s_{n-m_i}$, $s_0=1$,
$s_{n<0}=0$, where $p$, $v_i$, and $m_i$ are positive integers,
$m_1<\cdots<m_q$, and $r=0,\ldots,p-1$ \cite{AE-GenFibSqr}.  A
$(w,g;m)$-comb is a tile formed from a linear array of $m$ sub-tiles
(we call \textit{teeth}) each of size $w\times1$ (with the sides of
length $w$ aligned horizontally) and separated from one another by
gaps of size $g\times1$. A $(w,g;1)$-comb is just a $w\times1$
rectangular tile.  A $(w,g;2)$-comb is also known as a
$(w,g)$-fence. Fences were first introduced to provide a combinatorial
interpretation of the tribonacci numbers (which are defined by
$T_{n}=T_{n-1}+T_{n-2}+T_{n-3}+\delta_{n,2}$, $T_{n<2}=0$, where
$\delta_{i,j}$ is $1$ if $i=j$ and $0$ otherwise) via tiling an
$n$-board using just two types of tile, namely, squares and
$(\frac12,1)$-fences \cite{Edw08}.  Since then, we have used the
tiling of $n$-boards with various types of fence and rectangular or
square tiles to obtain combinatorial proofs of identities involving
the Fibonacci numbers \cite{EA19,EA20,EA20a,EA21}.

To obtain a combinatorial interpretation of the squares of various
number sequences via the tiling of an $n$-board with combs, we require
the following theorem which is a special case of Corollary~2.2 of
\cite{AE-GenFibSqr}.
\begin{thm}\label{T:genbij}
The number of ways to tile an $n$-board using
$(\frac12,\frac12;m_i)$-combs for $i=1,\ldots,q$ with
$0<m_1<\cdots<m_q$ is $s_n^2$ where
$s_n=s_{n-m_1}+\cdots+s_{n-m_q}+\delta_{n,0}$.
\end{thm}

For a combinatorial interpretation of the tribonacci numbers
squared we therefore need to consider tiling an $n$-board with
$(\frac12,\frac12;1)$-, $(\frac12,\frac12;2)$-, and
$(\frac12,\frac12;3)$-combs, which from now on we will refer to as
half-squares ($h$), fences ($f$), and combs ($c$), respectively. This
is formalized in the following theorem (which is a particular instance
of Theorem~\ref{T:genbij}).

\begin{thm}\label{T:bij}
Let $A_n$ be the number of ways to tile an $n$-board using
half-squares, fences, and combs. Then $A_n=T_{n+2}^2$.
\end{thm}

The tiling of a board using tiles with gaps (or tiles with any sides
which are not of integer length) can be treated as a tiling using
metatiles. A metatile is a grouping of tiles that exactly covers an
integer number of cells without any gaps and cannot be split into
smaller metatiles \cite{Edw08}. In some cases, there is a finite set
of possible metatiles; e.g., when tiling with squares and
$(\frac12,1)$-fences, the metatiles are a square, a fence with its gap
filled by a square, and three interlocking fences, and these are of
lengths 1, 2, and 3, respectively \cite{Edw08}.  Interesting
identities are obtained when there are metatiles of arbitrary length
if an expression for (or at least a recursion relation giving) the
number of metatiles of a given length can be obtained
\cite{EA19,EA20a,AE-GenFibSqr}. Obtaining such an expression has been
achieved via a digraph approach \cite{EA15}, determining the form of
all possible metatiles \cite{EA19}, finding a bijection between the
metatiles and some other objects whose number is known
\cite{AE-GenFibSqr}, and by examining how all possible metatiles of
length $l$ can be obtained by replacing tiles at the ends of metatiles
of length $l-1$ \cite{EA20a}. It is this final approach that we will
use in the next section. In Section~\ref{S:idn} we use the expressions
for the number of metatiles to obtain identities using similar methods
to \cite{EA20a}.

\section{Metatiles}\label{S:meta}

\begin{figure}[b]
\begin{center}
\includegraphics[width=11cm]{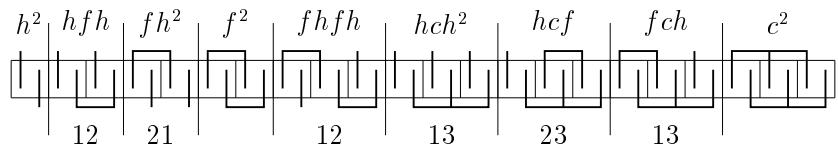}
\end{center}
\caption{A 22-board tiled with metatiles of length $l\leq3$. Each
  bold vertical line represents a ($\frac12\times1$) tooth and is
  aligned with the middle of the tooth it represents. Bold horizontal
  lines indicate which teeth are part of the same fence or comb.
  In the case of the $l=3$ mixed metatiles, only one element of
  each pair is shown. Thin vertical lines show boundaries between
  metatiles. The symbolic representation is above each
  metatile. $\sigma$ for the mixed metatiles appears below.}
\label{f:metatiles}
\end{figure}

As in \cite{EA20a,AE-GenFibSqr}, we define a \textit{mixed metatile}
as a metatile that contains more than one type of tile. Thus all the
metatiles are mixed except for $h^2$, $f^2$ (the bifence), and $c^2$
(the bicomb) which are shown in Fig.~\ref{f:metatiles}. We refer to
each half of each cell on the $n$-board as a
\textit{slot} \cite{EA20a,AE-GenFibSqr}.  The following lemma and its
proof are analogous to Lemma~2.3 in \cite{AE-GenFibSqr}.

\begin{lemma}\label{L:pair}
Each mixed metatile is a member of a pair of mixed metatiles. One
member of the pair is obtained from the other by swapping the contents
of the slots in each cell.
\end{lemma}
\begin{proof}
The swapping operation will only generate the same metatile if the
teeth in each pair of slots are from the same type of tile. This
can only occur if the metatile is not mixed. The swapping operation
will not split a metatile into more than one metatile since the
operation does not change which cells any given fence or comb
straddles.
\end{proof}

We use a 2-digit string $\sigma$ to specify the contents of the final
2 slots of a metatile. A digit $m$ in the string indicates that the
slot is filled by a tooth belonging to a
$(\frac12,\frac12;m)$-comb. Let $\mu_l^{[\sigma]}$ be the number of
metatiles of length $l$ that end with $\sigma$. Then by
Lemma~\ref{L:pair}, 
\begin{equation}\label{e:pair} 
\mu_l^{[m_1m_2]}=\mu_l^{[m_2m_1]},
\end{equation} 
where the $m_i$ are the digits of $\sigma$. As an example, $hfh$ and
$fh^2$ (Fig.~\ref{f:metatiles}) are a pair of mixed metatiles in the
sense of Lemma~\ref{L:pair}.

The digits in $\sigma$ for a mixed metatile must be distinct and, in
view of \eqref{e:pair}, we need only consider $\mu_l^{[12]}$,
$\mu_l^{[13]}$, and $\mu_l^{[23]}$. 

\begin{lemma}
For all integers $l$, when tiling with $h$, $f$, and $c$, 
\begin{subequations}
\label{e:m}
\begin{align} 
\label{e:m12} 
\mu_l^{[12]}&=\mu_{l-1}^{[12]}+\mu_{l-2}^{[12]}+\mu_{l-3}^{[12]}+\delta_{l,2}+\delta_{l,4},
\quad \mu_{l<2}^{[12]}=0, \\
\label{e:m13} 
\mu_l^{[13]}&=\mu_{l-1}^{[13]}+\mu_{l-2}^{[13]}+\mu_{l-3}^{[13]}+2\delta_{l,3},
\quad  \mu_{l<3}^{[13]}=0, \\
\label{e:m23} 
\mu_l^{[23]}&=\mu_{l-1}^{[23]}+\mu_{l-2}^{[23]}+\mu_{l-3}^{[23]}+\delta_{l,3}+\delta_{l,5},
\quad  \mu_{l<3}^{[23]}=0.
\end{align} 
\end{subequations}
\end{lemma}
\begin{proof}
To form metatiles of length $l$ from a given metatile of length $l-1$
we replace an $h$ in the final slot by an $f$, and/or replace an $f$
by a $c$. A remaining empty slot in the $l$th cell is filled with an
$h$. A metatile with $\sigma=12$ (e.g., $hfh$) will thus generate metatiles with
$\sigma=21$ (by replacing the final $h$ by an $f$ and adding an
$h$, e.g., $hf^2h$ which is paired with $fhfh$), $\sigma=13$ (by replacing the final $f$ by a $c$ and adding an
$h$ in the gap, e.g., $hch^2$), and $\sigma=23$ (by replacing both,
e.g., $hcf$). Similarly, metatiles with $\sigma=13$ and $\sigma=23$
will generate metatiles with $\sigma=21$ and $\sigma=31$,
respectively. Hence, after using \eqref{e:pair}, 
\begin{subequations}
\begin{align}
\label{e:pm12} 
\mu_l^{[12]}&=\mu_{l-1}^{[12]}+\mu_{l-1}^{[13]}+\delta_{l,2},\\
\label{e:pm13} 
\mu_l^{[13]}&=\mu_{l-1}^{[12]}+\mu_{l-1}^{[23]}+\delta_{l,3},\\
\label{e:pm23} 
\mu_l^{[23]}&=\mu_{l-1}^{[12]},
\end{align} 
\end{subequations}
where the $\delta_{l,2}$ and $\delta_{l,3}$ terms arise from the creation
of metatiles $hfh$ and $fch$ 
from the non-mixed metatiles $h^2$ and $f^2$, respectively.
Substituting \eqref{e:pm13} into \eqref{e:pm12} and then
\eqref{e:pm23} into the result gives \eqref{e:m12}. Then \eqref{e:m23}
immediately follows from applying \eqref{e:m12} to
\eqref{e:pm23}. Substituting \eqref{e:m12} and \eqref{e:m23} into
\eqref{e:pm13} and regrouping terms gives \eqref{e:m13}.
\end{proof}

Let $\mu_l$ be the number of mixed metatiles of length $l$. Then
$\mu_l=2(\mu_l^{[12]}+\mu_l^{[13]}+\mu_l^{[23]})$ and so from summing
\eqref{e:m} we obtain
\begin{equation}\label{e:mu_l}
\mu_l=\mu_{l-1}+\mu_{l-2}+\mu_{l-3}+6\delta_{l,3}
+2(\delta_{l,2}+\delta_{l,4}+\delta_{l,5}), \quad \mu_{l<2}=0.
\end{equation}
This gives $\{\mu_l\}_{l\ge2}=2,8,12,24,44,80,148,272,500,\ldots$,
which after a few terms is one form of tribonacci sequence. However,
$\mu_l$ can be expressed in terms of the usual tribonacci numbers as
shown in the following lemma.

\begin{lemma}\label{L:mu} 
The number of mixed metatiles of length $l$ when tiling with $h$, $f$,
and $c$ is
\[
\mu_l=4(T_{l}+T_{l-1})-2\delta_{l,2}.
\]
\end{lemma}
\begin{proof}
The generating function for \eqref{e:mu_l} is
$(2z^2+6z^3+2z^4+2z^5)/(1-z-z^2-z^3)$ which can be re-expressed as
$4(z^2+z^3)/(1-z-z^2-z^3)-2z^2$. As the generating function for $T_l$
is $z^2/(1-z-z^2-z^3)$, the result follows.
\end{proof}

\section{Identities}\label{S:idn}

\begin{lemma}\label{L:Anidn}
For all non-negative integers $n$,
\begin{equation}\label{e:Anidn}
A_{n}=\delta_{n,0}+A_{n-1}+3A_{n-2}+9A_{n-3}+\sum_{l=4}^n\mu_lA_{n-l},
\end{equation}
where $A_{n}=0$ for $n<0$.
\end{lemma}
\begin{proof}
Following \cite{BHS03,EA15,EA19}, we condition on the last
metatile. If the last metatile is of length $l$ there will be
$A_{n-l}$ ways to tile the remaining $n-l$ cells.  There is one
metatile of length 1 ($h^2$), three of length 2, nine of length 3,
and $\mu_l$ metatiles
of length $l$ for each $l\geq4$.  If $n=l$ there is exactly one tiling
corresponding to that final metatile so we make $A_0=1$. There is no
way to tile an $n$-board if $n<l$ and so $A_{n<0}=0$.
\end{proof}

\begin{idn}\label{I:sum}
For all non-negative integers $n$,
\[
T_{n}^2=\delta_{n,2}+T_{n-1}^2+3T_{n-2}^2+9T_{n-3}^2+4\sum_{l=4}^{n-2}(T_l+T_{l-1})T_{n-l}^2.
\]
\end{idn}
\begin{proof}
It follows from Lemma~\ref{L:Anidn}, Lemma~\ref{L:mu},
and Theorem~\ref{T:bij}.
\end{proof}

\begin{idn}\label{I:t2rr}
For $n\ge0$,
$T_{n}^2=\delta_{n,2}-\delta_{n,3}-\delta_{n,4}-\delta_{n,5}
+2T_{n-1}^2+3T_{n-2}^2+6T_{n-3}^2-T_{n-4}^2-T_{n-6}^2$.
\end{idn}
\begin{proof}
Letting $E(n)$ represent \eqref{e:Anidn} and re-indexing three of the
sums in $E(n)-E(n-1)-E(n-2)-E(n-3)$ gives
\begin{multline*} 
A_n=\delta_{n,0}-\delta_{n,1}-\delta_{n,2}-\delta_{n,3}
+2A_{n-1}+3A_{n-2}+6A_{n-3}-A_{n-4}-A_{n-6}\\
\mbox{}+\sum_{l=7}^n(\mu_l-\mu_{l-1}-\mu_{l-2}-\mu_{l-3})A_{n-l}.
\end{multline*} 
after using $\mu_4=12$, $\mu_5=24$, and $\mu_6=44$. The sum vanishes
by virtue of \eqref{e:mu_l} and, after changing $n$ to $n-2$,
the identity follows from Theorem~\ref{T:bij}.
\end{proof}

\begin{idn}\label{I:fc}
For all non-negative integers $n$,
\[
T_{n+4}^2=1+\sum_{k=0}^n\left\{3T_{k+2}^2+9T_{k+1}^2
+4\sum_{i=2}^k(T_{k+4-i}+T_{k+3-i})T_i^2\right\}.
\]
\end{idn}
\begin{proof}
How many ways are there to tile an $(n+2)$-board using at least one
fence or comb? \textit{Answer~1}: $A_{n+2}-1$ since this corresponds
to all tilings except the all-$h$ tiling.  \textit{Answer~2}:
condition on the location of the last metatile containing an $f$ or
$c$. If the metatile is of length $l$ and ends on cell $k+2$ (for
$k=l-2,\ldots,n$) then the number of ways to tile the board is 
$(\mu_l+\delta_{l,2}+\delta_{l,3})A_{k+2-l}$. Summing over all
possible lengths and all possible $k$, introducing $i=k+2-l$ and
summing over that rather than $l$, and then equating to Answer~1 gives
\[
A_{n+2}-1=\sum_{k=0}^n\sum_{i=0}^k(\mu_{k+2-i}+\delta_{i,k}+\delta_{i,k-1})A_i.
\]
Taking the $i=k$ and $i=k-1$ terms outside the sum over $i$, changing
$i$ to $i-2$, and then using Theorem~\ref{T:bij} and Lemma~\ref{L:mu}
gives the identity.
\end{proof}

\begin{idn}\label{I:hc}
For all non-negative integers $n$ and $j=0,1$,  
\[
T_{2(n+1)+j}^2=1+\sum_{k=1}^n\left\{T_{2k+j+1}^2+2T_{2k+j}^2+9T_{2k+j-1}^2
+4\!\!\!\sum_{i=0}^{2k+j-4}(T_{2k+j-i}+T_{2k+j-i-1})T_{i+2}^2\right\}.
\]
\end{idn}
\begin{proof}
How many ways are there to tile an $(2n+j)$-board using at least one
half-square or comb? \textit{Answer~1}: $A_{2n+j}-\delta_{0,j}$ since
only the all-bifence tiling has no $h$ or $c$ and this only occurs for
even length boards. \textit{Answer~2}:
condition on the location of the last metatile containing an $h$ or
$c$. The last cell of this metatile must lie on cell $2k+j$ for some
$k=\delta_{0,j},\ldots,n$ since the cells to the right must be filled
with bifences. If the metatile is of length $l$
 then the number of ways to tile the board is 
$(\mu_l+\delta_{l,1}+\delta_{l,3})A_{2k+j-l}$. Summing over all
possible $l$ and $k$, separating out the $l=1$ case, and then 
equating to Answer~1 gives
\[
A_{2n+j}-\delta_{j,0}=\sum_{k=\delta_{j,0}}^nA_{2k+j-1}
+\sum_{k=1}^n\sum_{i=0}^{2k+j-2}(\mu_{2k+j-i}+\delta_{3,2k+j-i})A_i.
\]
Simplifying and then using Theorem~\ref{T:bij} and Lemma~\ref{L:mu} gives
the identity.
\end{proof}

\begin{idn}\label{I:hf}
For all non-negative integers $n$ and $j=0,1,2$,  
\[
T_{3n+2+j}^2=1+3\delta_{j,2}+\sum_{k=1}^n\left\{T_{3k+j+1}^2+3T_{2k+j}^2
+4\!\!\!\sum_{i=0}^{3k+j-3}(T_{3k+j-i}+T_{3k+j-i-1})T_{i+2}^2\right\}.
\]
\end{idn}
\begin{proof}
How many ways are there to tile an $(2n+j)$-board using at least one
half-square or fence? \textit{Answer~1}: $A_{3n+j}-\delta_{0,j}$ since
only the all-bicomb tiling has no $h$ or $f$ and this only occurs for
boards of length divisible by 3. \textit{Answer~2}: condition on the
location of the last metatile containing an $h$ or $f$. The last cell
of this metatile must lie on cell $3k+j$ for some
$k=\delta_{0,j},\ldots,n$ since the cells to the right must be filled
with bicombs. If the metatile is of length $l$ then the number of ways
of tile the board is
$(\mu_l+\delta_{l,1}+\delta_{l,2})A_{3k+j-l}$. Summing over all
possible $l$ and $k$, separating out the $l=1$ and $l=2$ cases,
equating the whole expression to
Answer~1, simplifying, and then using Theorem~\ref{T:bij} and
Lemma~\ref{L:mu} gives the identity.
\end{proof}

The number of ways to tile an $n$-board using only $h^2$, $f^2$, and
$c^2$ is $T_{n+2}$ since these metatiles are of lengths 1, 2, and 3.

\begin{idn}\label{I:mm}
For all non-negative integers $n$,
\[
T_n^2=T_n+\sum_{k=2}^{n-2}\sum_{l=2}^k\{4(T_l+T_{l-1})-2\delta_{l,2}\}
T_{k-l+2}^2T_{n-k}.
\]
\end{idn}
\begin{proof}
How many ways are there to tile an $(n-2)$-board using at least one
mixed metatile?  \textit{Answer~1}: $A_{n-2}-T_n$ since $T_n$ is the
number of ways to tile an $(n-2)$-board without using mixed
metatiles. \textit{Answer~2}: condition on the last mixed metatile. If
this metatile is of length $l$ and ends on cell $k$ where
$k=l,\ldots,n-2$, then the number of ways to tile the remaining cells
is $A_{k-l}T_{n-2-k}$. Given that there are $\mu_l$ mixed metatiles of
length $l$, summing over all possible $k$ and $l$ and equating to
Answer~1 gives
\[
A_{n-2}-T_n=\sum_{k=2}^{n-2}\sum_{l=2}^k\mu_lA_{k-l}T_{n-2-k}.
\]
Using Theorem~\ref{T:bij} and Lemma~\ref{L:mu} gives the identity.
\end{proof}

The number of ways to tile an $n$-board using only $h$ and $f$ is
$F_{n+1}^2$ and there are 2 metatiles of length $l$ containing only
$h$ and $c$ for $l\ge3$ \cite{EA19}.
\begin{idn}\label{I:c}
For all non-negative integers $n$,
\[
T_{n+2}^2=F_{n+1}^2+\sum_{k=3}^n\sum_{l=3}^k\{4(T_l+T_{l-1})+\delta_{l,3}-2\}
T_{k-l+2}^2F_{n-k+1}^2.
\]
\end{idn}
\begin{proof}
How many ways are there to tile an $n$-board using at least one comb?
\textit{Answer~1}: $A_n-F_{n+1}^2$. \textit{Answer~2}: condition on
the last metatile containing a comb. The number of such metatiles of
length $l$ (where, owing to the size of a comb, $l\ge3$) is
$\mu_l+\delta_{l,3}-2$. Suppose such a metatile of length $l$ ends on
cell $k$ (where $k=l,\ldots,n$). Then the number of ways to tile the
rest of the board is $A_{k-l}F_{n-k+1}^2$. Summing over all possible
metatiles, $l$, and $k$, and equating to Answer~1 gives
\[
A_n-F_{n+1}^2=\sum_{k=3}^{n}\sum_{l=3}^k(\mu_l+\delta_{l,3}-2)A_{k-l}F^2_{n-k+1}.
\]
Using Theorem~\ref{T:bij} and Lemma~\ref{L:mu} gives the identity.
\end{proof}

The Narayana's cows sequence is defined by
$c_n=c_{n-1}+c_{n-3}+\delta_{n,0}$, $c_{n<0}=0$.  From
Theorem~\ref{T:genbij}, the number of ways to tile an $n$-board using
only $h$ and $c$ is $c_n^2$. The number of mixed metatiles of length $l$
containing only $h$ and $c$ is $2p_{l-1}$
where $p_n=p_{n-2}+p_{n-3}+\delta_{n,0}$, $p_{n<0}=0$
\cite{AE-GenFibSqr}. These are the Padovan numbers (offset from
sequence A000931 in \cite{Slo-OEIS}). From Theorem~\ref{T:genbij}, the
number of ways to tile an $n$-board using only $f$ and $c$ is
$p_n^2$. The number of mixed metatiles of length $l$ containing only $f$ and
$c$ is $2c_{l-5}$ \cite{AE-GenFibSqr}.

\begin{idn}\label{I:f}
For all non-negative integers $n$,
\[
T_{n+2}^2=c_n^2+\sum_{k=2}^n\sum_{l=2}^k\{4(T_l+T_{l-1})-\delta_{l,2}-2p_{l-1}\}
T_{k-l+2}^2c_{n-k}^2.
\]
\end{idn}
\begin{proof}
How many ways are there to tile an $n$-board using at least one fence?
\textit{Answer~1}: $A_n-c_n^2$. \textit{Answer~2}: condition on the
last metatile containing a fence. The number of such metatiles of
length $l$ is $\mu_l+\delta_{l,2}-2p_{l-1}$. The proof then proceeds
in an analogous way to that of Identity~\ref{I:c}.
\end{proof}

\begin{idn}\label{I:h}
For all non-negative integers $n$,
\[
T_{n+2}^2=p_n^2+\sum_{k=1}^n\sum_{l=1}^k\{4(T_l+T_{l-1})+\delta_{l,1}-2\delta_{l,2}-2c_{l-5}\}
T_{k-l+2}^2p_{n-k}^2.
\]
\end{idn}
\begin{proof}
How many ways are there to tile an $n$-board using at least one
half-square?  \textit{Answer~1}: $A_n-p_n^2$. \textit{Answer~2}:
condition on the last metatile containing a half-square. The number of
such metatiles of length $l$ is $\mu_l+\delta_{l,1}-2c_{l-5}$. The
proof then proceeds as in Identity~\ref{I:c}.
\end{proof}

\begin{lemma}\label{L:m1m2}
The number of ways to tile a board using $h$, $f$, and $c$ such that
$\sigma=m_1m_2$ for $m_1,m_2\in\{1,2,3\}$ is $T_{n+2-m_1}T_{n+2-m_2}$.
\end{lemma}
\begin{proof}
There is a bijection between the tiling of an ordered pair of
$n$-boards with squares, dominoes, and trominoes and the tiling of an
$n$-board with $h$, $f$, and $c$ (see the proof of Theorem~2.1 in
\cite{AE-GenFibSqr}) whereby all teeth of a
$(\frac12,\frac12;m)$-comb ending at the left (right) slot the $k$th
cell correspond to an $m$-omino ending in the $k$th cell of the first
(second) of the pair of boards. Thus when there is a tooth
belonging to a $(\frac12,\frac12;m_i)$-comb in the final left (right)
slot there are $T_{n-m_i+2}$ ways to tile the rest of the first
(second) board ending in an $m_i$-omino and hence
$T_{n+2-m_1}T_{n+2-m_2}$ ways to tile the remaining cells of both
boards with ominoes.
\end{proof}

\begin{idn}\label{I:TnT} For integers $n\ge2$, 
\[
T_{n+1}T_n=\sum_{l=2}^n(T_l+T_{l-2})T_{n-l+2}^2. 
\]
\end{idn}
\begin{proof}
How many ways are there to tile an $n$-board with an $h$ in the final
left slot and an $f$ tooth in the final right slot? \textit{Answer 1}: by
Lemma~\ref{L:m1m2}, $T_{n+1}T_n$. \textit{Answer 2}: there are
$\mu_l^{[12]}$ metatiles of length $l$ than can end such a board and
thus, summing over all possible $l$, the number of ways to tile the board is
$\sum_{l=2}^n\mu_l^{[12]}A_{n-l}$. Equating this to Answer~1 and using
Theorem~\ref{T:bij} and the
result from \eqref{e:m12} that $\mu_l^{[12]}=T_l+T_{l-2}$ then gives
the identity.
\end{proof}

\begin{idn} For integers $n\ge3$, 
\[
T_{n+1}T_{n-1}=2\sum_{l=3}^nT_{l-1}T_{n-l+2}^2. 
\]
\end{idn}
\begin{proof}
How many ways are there to tile an $n$-board with an $h$ in the final
left slot and a $c$ tooth in the final right slot? \textit{Answer 1}: by
Lemma~\ref{L:m1m2}, $T_{n+1}T_{n-1}$. \textit{Answer 2}: there are
$\mu_l^{[13]}$ metatiles of length $l$ than can end such a board. 
From \eqref{e:m13}, $\mu_l^{[13]}=2T_l$. The proof then proceeds
in an analogous way to that for Identity~\ref{I:TnT}.
\end{proof}

\section{Discussion}
The unusual feature of this particular selection of tile types is that
numbers of metatiles can be simply expressed in terms of the sequence
generated from the tiling of the whole board. As a result, most of the
identities we obtain involve only the tribonacci numbers.  Aside from
Identity~\ref{I:t2rr}, we have not been able to locate any of the
remaining identities we derive here in the literature.


\medskip

\noindent MSC2010:  05A19, 11B39

\end{document}